\documentclass[a4paper,10pt,reqno]{amsart}
\usepackage{amsmath,amsfonts,amsthm}
\usepackage[mathscr]{eucal}
\usepackage{rotating}
\usepackage{multirow}
\usepackage{slashbox}

\numberwithin{equation}{section}
\newtheorem{theorem}{Theorem}[section]

\newcommand{\RR}{\mathbb{R}}

\begin{document}

\title[Solutions of SDEs obeying LIL]{Recurrent Solutions of Stochastic Differential Equations with Non-constant Diffusion Coefficients which obey the Law of the Iterated Logarithm} 

\author{John A. D. Appleby}
\address{Edgeworth Centre for Financial Mathematics, School of Mathematical Sciences, Dublin City University,
Ireland} \email{john.appleby@dcu.ie} \urladdr{webpages.dcu.ie/\textasciitilde applebyj}

\author{Huizhong Appleby-Wu}
\address{Department of Mathematics, St Patrick's College, Drumcondra, Dublin 9, Ireland} \email{huizhong.applebywu@spd.dcu.ie}


\thanks{JA thanks SFI for the support of this research under the Mathematics
Initiative 2007 grant 07/MI/008 ``Edgeworth Centre for Financial
Mathematics''.}



\date{27 October 2012}

\begin{abstract}
By using a change of scale and space, we study a class of stochastic
differential equations (SDEs) whose solutions are drift--perturbed
and exhibit behaviour analogous to standard Brownian motion
including to the Law of the Iterated Logarithm (LIL). Sufficient
conditions ensuring that these processes obey the LIL are given.
\end{abstract}

\maketitle

\section{Introduction}

The Law of the Iterated Logarithm is one of the most
important results on the asymptotic behavior of one--dimensional
standard Brownian motion:
\begin{equation}\label{LILBM}
\limsup_{t\to\infty}\frac{|B(t)|}{\sqrt{2t\log\log{t}}}=1,\quad\text{a.s.}
\end{equation}
Classical work on iterated logarithm--type results, as well as
associated lower bounds on the growth of transient processes date
back to Dvoretzky and Erd\H{o}s~\cite{DvErd}. There is an
interesting literature on iterated logarithm results and the growth
of lower envelopes for self-similar Markov processes (cf. e.g.,
Rivero~\cite{Riv}, Chaumont and Pardo~\cite{ChPardo}) which exploit
a Lamperti representation~\cite{Lamp}, processes conditioned to
remain positive (cf. Hambly et al.~\cite{Hamblyetal}), and diffusion
processes with special structure (cf. e.g. Bass and
Kumagi~\cite{BasKum}).

In contrast to these papers the analysis here is inspired by work
of Motoo~\cite{Motoo} on iterated logarithm results for Brownian
motions in finite dimensions, in which the asymptotic behaviour is
determined by means of time change arguments to reduce the process
under study to a stationary one.

\section{Main Result}

\subsection{Preliminaries}
Throughout the paper, we use
$(\Omega,\mathcal{F},\{\mathcal{F}(t)\}_{t\geq 0},\mathbb{P})$ to
denote a complete filtered probability space. The set of
non-negative real numbers is denoted by $\RR^+$. Let $L^1([a,b];\RR)$ be the family of Borel
measurable functions $h:[a,b]\to \RR$ such that
$\int_a^b|h(x)|dx<\infty$. The abbreviation \emph{a.s.} stands for
\emph{almost surely}. 

Throughout the paper, we assume that both the drift and the
diffusion coefficients of SDE being studied satisfies the local
Lipschitz condition. If an
autonomous scalar SDE has drift coefficient $f(\cdot)$ and
non-degenerate diffusion coefficient $g(\cdot)$, then a scale
function and speed measure of solution of this SDE are defined as
\begin{equation}\label{scale}
s_c(x)=\int_c^xe^{-2\int_c^y\frac{f(z)}{g^2(z)}\,dz}\,dy,\quad m(dx)=\frac{2dx}{s'(x)g^2(x)},\quad c,x\in (0,\infty)
\end{equation}
respectively. These functions help us to determine the recurrence
and stationarity of a process on $(0,\infty)$ (cf.\cite{K&S}). Moreover, Feller's test for explosions (cf.\cite{K&S})
allows us to examine whether a process will never escape from its
state space in finite time. This in turn relies on whether
\[
v(0+)=v(\infty-)=\infty
\]
or not, where $v$ is defined as
\begin{equation}\label{vdef}
v_c(x)=\int_c^xs'_c(y)\int_c^y\frac{2dz}{s'_c(z)g^2(z)}\,dy,\quad c\in (0,\infty),\quad x\in(0,\infty).
\end{equation}

As mentioned in the introduction, Motoo's Theorem is an important
tool in determining the largest deviations for any stationary or
asymptotically stationary processes, we state it here for future
use.
\begin{theorem}~\label{thm.motoo}[Motoo]
Let $X$ be the unique continuous real valued process satisfying the
following equation
\begin{gather*}
dX(t)=f(X(t))\,dt+g(X(t))\,dB(t),\quad t\geq 0,
\end{gather*}
with $X(0)=x_0$. Let $s$ and $m$ be the scale function and speed
measure of $X$ as defined in \eqref{scale}, and let $h:(0,\infty)\to
(0,\infty)$ be an increasing function. If $X$ is recurrent on $(0,\infty)$ (or $[0,\infty)$ in the case when $0$ is an instantaneous reflecting point)
and $m(0,\infty)<\infty$, then
\begin{equation*}
\mathbb{P}\left[\limsup_{t\to\infty}\frac{X(t)}{h(t)}\geq
1\right]=1\,\,\text{or}\,\,0
\end{equation*}
according to whether
\begin{equation*}
\int_{t_0}^{\infty}\frac{1}{s(h(t))}\,dt=\infty\quad\text{or}\quad \int_{t_0}^{\infty}\frac{1}{s(h(t))}\,dt<\infty,\quad\text{for some}\,\,t_0>0.
\end{equation*}
\end{theorem}

\subsection{The main result}
Let $g:\mathbb{R}\to\mathbb{R}$ be locally Lipschitz continuous and suppose that $f:[0,\infty)\times \mathbb{R}\to\mathbb{R}$ 
is locally Lipschitz continuous. Suppose that 
\begin{gather}\label{xfxbound}
\text{There exists $\rho>0$ such that }
\sup_{(x,t)\in \mathbb{R}\times [0,\infty)}xf(x,t)\leq\rho,\\
 \label{xfxlowerbound}
 \text{There exists $\mu>-1/2$ such that }
\mu:= \inf_{(x,t)\in \mathbb{R}\times [0,\infty)}\frac{xf(x,t)}{g^2(x)}.
 \end{gather}
Suppose also that $g:\mathbb{R}\to\mathbb{R}$ satisfies
\begin{equation}\label{co5.4g}
\forall\,x\in \mathbb{R},\quad g(x)\neq 0,\quad \lim_{|x|\to\infty}g(x)=\sigma \in\mathbb{R}/\{0\}.
\end{equation}
Then there exists a unique continuous adapted process $X$ which obeys 
\begin{equation}\label{eq.sde}
dX(t)=f(X(t),t)\,dt+g(X(t))\,dB(t),\quad t\geq 0, \quad X(0)=x_0
\end{equation}
(see e.g., \cite{Mao}). 

\begin{theorem} \label{thm.mainresult}
Suppose that $f:[0,\infty)\times \mathbb{R}\to\mathbb{R}$ is locally Lipschitz continuous.
Suppose $g:\mathbb{R}\to\mathbb{R}$ is locally Lipschitz, even and satisfies \eqref{co5.4g}, and 
that $f$ and $g$ obey \eqref{xfxbound}. Then there is a unique continuous adapted process 
satisfying \eqref{eq.sde}. If moreover $f$ and $g$ obey \eqref{xfxlowerbound}, then $X$ obeys 
 \begin{equation} \label{bound}
 \limsup_{t\to\infty} \frac{|X(t)|}{\sqrt{2t\log\log t}} =|\sigma|, \quad \text{a.s.}
 \end{equation}
\end{theorem}

In a recent paper~\cite{AppWu:2008b}, we established \eqref{bound} for the solution of \eqref{eq.sde} under the conditions \eqref{xfxbound}, \eqref{co5.4g}.
In the condition \eqref{xfxlowerbound} it was presumed that $\mu>1/2$. This restriction forces $|X(t)|\to\infty$ as $t\to\infty$ a.s. 
Therefore we have shown here that the asymptotic growth implied by \eqref{bound} also holds in some cases where the 
diffusion coefficient is asymptotically constant and the solution can be \emph{recurrent} rather than \emph{transient}. 

In the case when $g$ is constant with $g(x)=\sigma$ (so that $g$ obviously obeys \eqref{co5.4g}), $f$ obeys \eqref{xfxbound} and \eqref{xfxlowerbound} (with the last condition written as $\inf_{(x,t)\in\mathbb{R}\times[0,\infty)} xf(x,t)>-\sigma^2/2$), 
it was proven in~\cite{AppWu:2008b} that the solution of \eqref{eq.sde} obeys \eqref{bound}. In this paper, we have managed to extend the result to 
the case where the diffusion coefficient is only \emph{asymptotically} constant. 

\section{Proof of Theorem~\ref{thm.mainresult}}
The result is proven by establishing that the solution $X$ obeys 
\begin{equation} \label{upbound}
 \limsup_{t\to\infty} \frac{|X(t)|}{\sqrt{2t\log\log t}}\leq |\sigma|, \quad \text{a.s.}
 \end{equation}
 and then that it obeys 
 \begin{equation} \label{lowbound}
 \limsup_{t\to\infty} \frac{|X(t)|}{\sqrt{2t\log\log t}}\geq |\sigma|, \quad \text{a.s.}
 \end{equation}
 The proof of \eqref{upbound} was given in \cite{AppWu:2008b}, but is included here in order to 
 introduce notation used in the proof of \eqref{lowbound}.  

\subsection{Proof of \eqref{upbound}}
Without loss of generality, we can choose $\rho>\sigma^2/2$. 
By It\^o's rule 
\[
X^2(t)=x_0^2+\int_0^t 2X(s)f(X(s),s)+g^2(X(s))\,ds+\int_0^t 2X(s)g(X(s))\,dB(s),\quad\text{a.s.}
\]
Let $Z(t)=X^2(t)$ so that $|X(t)|=\sqrt{Z(t)}$ for all $t\geq 0$. Define 
$\gamma(x)=x/|x|$ for $x\neq 0$ and $\gamma(0)=1$. Then $\gamma^2(x)=1$ for all $x\in\mathbb{R}$. 
Also define 
\[
W(t)=\int_0^t \gamma(X(s))\,dB(s),\quad\text{a.s.}
\]
Then $W$ is a $\mathcal{F}^B$--adapted Brownian motion. If $Z(t)\neq 0$ 
\begin{align*}
2\sqrt{Z(t)}g(\sqrt{Z(t)})\gamma(X(t))&=2|X(t)|g(|X(t)|)\gamma(X(t))=2|X(t)|g(X(t))\gamma(X(t))\\
&=2X(t)g(X(t))\gamma^2(X(t))=2X(t)g(X(t)).
\end{align*}
If $Z(t)=0$ then $2\sqrt{Z(t)}g(\sqrt{Z(t)})\gamma(X(t))=0=2X(t)g(X(t))$. 
Hence,  
\begin{equation} \label{eq.z}
Z(t)=x_0^2+\int_0^t 2X(s)f(X(s),s) +g^2(\sqrt{Z(s)})\,ds+\int_0^t 2\sqrt{Z(s)}g(\sqrt{Z(s)})\,dW(s), 
\end{equation}
and so with $\beta(t):=2(X(t)f(X(t),t)-\rho)\leq 0$ for $t\geq 0$, we have 
\[
Z(t)=x_0^2+\int_0^t \beta(s) + 2\rho +g^2(\sqrt{Z(s)})\,ds+\int_0^t 2\sqrt{Z(s)}g(\sqrt{Z(s)})\,dW(s), \quad\text{a.s.}
\]
Note also that $\beta$ is $\mathcal{F}^B$-adapted. Define $Z_u(t)$ by 
\begin{equation} \label{eq.zu}
Z_u(t)=1+x_0^2+\int_0^t 2\rho +g^2(\sqrt{|Z_u(s)|})\,ds+\int_0^t 2\sqrt{|Z_u(s)|}g(\sqrt{|Z_u(s)|})\,dW(s), \quad\text{a.s.}
\end{equation}
Since $W$ is a $\mathcal{F}^B$--Brownian motion, $Z_u$ is $\mathcal{F}^B$--adapted. Notice that 
\[
dZ_u(t)=\left\{2\rho +g^2(\sqrt{|Z_u(t)|})\right\}\,dt+2\sqrt{|Z_u(t)|}g(\sqrt{|Z_u(t)|})\,dW(t).
\]
Define the process $X_u$ by  
\begin{equation*}
dX_u(t)=\frac{\rho}{X_u(t)}\,dt+g(X_u(t))\,dB(t),\quad t\geq 0
\end{equation*}
where $X_u(0)=\sqrt{1+|x_0|^2}$. It is easy to check that the scale function of $X_u$ satisfies $s_{X_u}(\infty)<\infty$ and 
$s_{X_u}(0)=-\infty$. Thus $\mathbb{P}\,[\lim_{t\to\infty}X_u(t)=\infty]=1$. Moreover $v_{X_u}(\infty)=v_{X_u}(0)=\infty$, 
which implies that $\mathbb{P}\,[X_u(t)>0;\,\forall\,0<t<\infty]=1$. Then $Z_u=X_u^2>0$ obeys \eqref{eq.zu}. Therefore 
we have 
\[
dZ_u(t)=\left\{2\rho +g^2(\sqrt{Z_u(t)})\right\}\,dt+2\sqrt{Z_u(t)}g(\sqrt{Z_u(t)})\,dW(t),
\]
and so by the Ikeda--Watanabe comparison theorem~\cite[Chapter VI, Theorem 1.1]{IkWat} we have 
$Z_u(t)\geq Z(t)$ for all $t\geq 0$ a.s. and so $X^2(t)\leq X_u^2(t)$ for all $t\geq 0$ a.s.
By a result in \cite{AppWu:2008b}, it is known that $X_u$ obeys
\begin{equation*}
\limsup_{t\to\infty}\frac{X_u(t)}{\sqrt{2t\log\log{t}}}=|\sigma|,\quad \text{a.s.},
\end{equation*}
which implies \eqref{upbound}. 

\subsection{Proof of \eqref{lowbound}}
Define 
\[
\beta_2(t)=2\left(\frac{X(t)f(X(t),t)}{g^2(X(t))} - \mu\right)g^2(X(t)), \quad t \geq 0.
\] 
Then $\beta_2$ is $\mathcal{F}^B$--adapted and $\beta_2(t)\geq 0$ for all $t\geq 0$. Moreover 
\[
2X(t)f(X(t),t)+g^2(X(t))= \beta_2(t) +(2\mu +  1) g^2(\sqrt{Z(t)}), \quad t\geq 0.
\]
Therefore from \eqref{eq.z} we obtain
\begin{equation} \label{eq.zzl}
Z(t)=x_0^2+\int_0^t \left\{\beta_2(s) + (2\mu+  1) g^2(\sqrt{Z(s)})\right\}\,ds+\int_0^t 2\sqrt{Z(s)}g(\sqrt{Z(s)})\,dW(s), 
\end{equation}
Define $Z_L(t)$ by 
\begin{equation} \label{eq.zl}
Z_L(t)=x_0^2+\int_0^t (2\mu+  1) g^2(\sqrt{|Z_L(s)|})\,ds+\int_0^t 2\sqrt{|Z_L(s)|}g(|\sqrt{Z_L(s)|})\,dW(s), \quad\text{a.s.}
\end{equation}
Define 
\begin{equation} \label{eq.theta}
\theta(t)=\int_0^t g^2(\sqrt{|Z_L(s)|})\,ds, \quad t\geq 0.
\end{equation}
Then $\theta$ is increasing and let $\tau=\theta^{-1}$. Define $\tilde{Z}_L(t)=Z_L(\tau(t))$ for $t\geq 0$. Let 
$\mathcal{G}(t)=\mathcal{F}^B(\tau(t))$ for all $t\geq 0$ and define 
\[
M(t)=\int_0^{\tau(t)} 2\sqrt{|Z_L(s)|}g(|\sqrt{Z_L(s)|})\,dW(s), \quad t\geq 0. 
\]
Then $M$ is a $\mathcal{G}(t)$--martingale. Clearly 
\[
\langle M\rangle(t)=\int_0^{\tau(t)} 4|Z_L(s)|g^2(|\sqrt{Z_L(s)|})\,ds = \int_0^t 4|\tilde{Z}_L(u)|\,du.
\]
Thus, by Theorem 3.4.2 in \cite{K&S}, there is an extension of  $(\tilde{\Omega},\tilde{\mathcal{F}},\tilde{\mathbb{P}})$ of $(\Omega,\mathcal{F},\mathbb{P})$ on which is defined a one--dimensional Brownian motion $\tilde{W}=\{\tilde{W}(t);\tilde{\mathcal{G}}(t);0\leq t<+\infty\}$ such that
\begin{equation} \label{eq.martingalezl}
M(t) = \int_0^t  2\sqrt{|\tilde{Z}_L(s)|} \,d\tilde{W}(s),
\quad\text{$\widetilde{\mathbb{P}}$-a.s.}
\end{equation}
The filtration $\tilde{\mathcal{G}}(t)$ in the extended space is such that $\tilde{Z}_L$ is $\tilde{\mathcal{G}}(t)$--adapted.  
Thus by \eqref{eq.zl} and \eqref{eq.martingalezl} we have
\begin{equation} \label{eq.zltilde}
\tilde{Z}_L(t) =x_0^2 + \int_0^t (2\mu+1)\,ds + \int_0^t  2\sqrt{|\tilde{Z}_L(s)|} \,d\tilde{W}(s),\quad t\geq 0.
\end{equation}
Since $\mu>-1/2$, it now follows that $\tilde{Z}_L(t)\geq 0$ for all $t\geq 0$ a.s. Therefore $Z_L(t)\geq 0$ for all $t\geq 0$ a.s. 
Hence $Z_L$ obeys 
\begin{equation*} 
dZ_L(t)= (2\mu+  1) g^2(\sqrt{Z_L(t)})\,dt+ 2\sqrt{Z_L(t)}g(\sqrt{Z_L(t)})\,dW(t),  \quad\text{$t\geq 0$ a.s.}
\end{equation*}
and so by \eqref{eq.zzl} and the Ikeda--Watanabe comparison theorem~\cite[Chapter VI, Theorem 1.1]{IkWat}, we have $Z(t)\geq Z_L(t)$ for all $t\geq 0$ a.s. Hence $X^2(\tau(t))\geq Z_L(\tau(t))=\tilde{Z}_L(t)$
for all $t\geq 0$ a.s. If we define $\bar{Z}_L(t)=\tilde{Z}_L(e^t-1)$ for $t\geq 0$, then there is another Brownian motion $\bar{W}$ such that 
\[
d\bar{Z}_L(t)=(2\mu+1)e^t \,dt + 2\sqrt{\bar{Z}_L(t)}e^{t/2}\,d\bar{W}(t), \quad t\geq 0.
\] 
Define 
\begin{equation} \label{eq.u}
U(t)=e^{-t}\tilde{Z}_L(e^t-1)=e^{-t}\bar{Z}_L(t), \quad t\geq 0. 
\end{equation}
Then $U(0)=x_0^2$ and 
\[
dU(t)=\left( (2\mu +1) - U(t) \right)\,dt + 2\sqrt{U(t)}\,d\bar{W}(t), \quad t\geq 0.
\]   
By theorem~\ref{thm.motoo}, it follows that 
\[
\limsup_{t\to\infty} \frac{U(t)}{2\log t} =1, \quad \text{a.s.}
\]   
Using the connection between $U$ and $\tilde{Z}_L$ we obtain
\[
\limsup_{t\to\infty} \frac{\tilde{Z}_L(t)}{2t\log\log t}=1, \quad \text{a.s.}
\]
Therefore
\[
\limsup_{t\to\infty} \frac{X^2(\tau(t))}{2t\log\log t}\geq \limsup_{t\to\infty} 
\frac{\tilde{Z}_L(t)}{2t\log\log t}=1, \quad \text{a.s.}
\]
Since $\theta=\tau^{-1}$ and $\theta(t)\to\infty$ as $t\to\infty$ we have  
\begin{equation} \label{eq.xsqlower1}
\limsup_{t\to\infty} \frac{X^2(t)}{2\theta(t)\log\log \theta(t)}\geq 1, \quad \text{a.s.}
\end{equation}
By \eqref{co5.4g}, $g^2(x)\to \sigma^2>0$ as $|x|\to\infty$. Since $g^2$ is continuous and $g(x)\neq 0$ for all $x\in\mathbb{R}$ (by assumption
\eqref{co5.4g}), it follows that there exist $K_1^2>0$ and $K_2^2\in[K_1^2,\infty)$ such that 
$0<K_1^2\leq g^2(x)\leq K_2^2$ for all $x\in\mathbb{R}$. Therefore
\begin{equation} \label{eq.thetabounds}
0<K_1^2 t \leq \theta(t)\leq K_2^2 t, \quad \text{for all $t\geq 0$ a.s.},
\end{equation} 
which implies 
\[
\lim_{t\to\infty} \frac{\log \log \theta(t)}{\log \log t}=1, \quad\text{a.s.}
\]
Using \eqref{eq.xsqlower1} now yields
\begin{equation} \label{eq.xsqlower2}
\limsup_{t\to\infty} \frac{X^2(t)}{2\theta(t)\log\log t}\geq 1, \quad \text{a.s.}
\end{equation}
The rest of the proof is devoted to showing that 
\begin{equation} \label{eq.thetaaverage}
\lim_{t\to\infty} \frac{1}{t}\int_0^t g^2(\sqrt{Z_L(s)})\,ds = \sigma^2,\quad\text{a.s.}
\end{equation}
which together with \eqref{eq.xsqlower2} yields \eqref{lowbound}.

We prove \eqref{eq.thetaaverage} in four steps:
\begin{itemize}
\item[{\bf Step A.}] If $U$ is the process defined in \eqref{eq.u}, then for every $c>0$ 
\begin{equation} \label{eq.uprobint}
\int_1^\infty \mathbb{\tilde{P}}[U(s)\leq ce^{-s}]\,ds < +\infty
\end{equation} 
\item[{\bf Step B.}] If $U$ is the process defined in \eqref{eq.u}, then \eqref{eq.uprobint} implies 
\begin{equation} \label{eq.conv1}
\lim_{t\to\infty} \frac{1}{e^t} \int_1^t e^s I_{\{U(s)\leq ce^{-s}\}}\,ds = 0, \quad \text{a.s.}
\end{equation}
\item[{\bf Step C.}] \eqref{eq.conv1} implies 
\begin{equation} \label{eq.conv3}
\lim_{t\to\infty} \frac{1}{t} \int_1^t I_{\{Z_L(s)\leq c\}}\,ds = 0, \quad \text{a.s.}
\end{equation}
\item[{\bf Step D.}] If $Z_L$ obeys \eqref{eq.conv3}, then it also obeys \eqref{eq.thetaaverage}. 
\end{itemize}
The proof of Steps A--D are given in the next four subsections, which completes the proof of \eqref{lowbound}.


\subsection{Proof of Step A, i.e., \eqref{eq.uprobint}}
Clearly by the definition of $U$ and $\tilde{Z}_L$ we have 
\begin{equation} \label{eq.uprobzlprob}
\int_1^\infty \mathbb{\tilde{P}}[U(s)\leq ce^{-s}]\,ds = \int_{e-1}^\infty \frac{1}{1+t}\mathbb{\tilde{P}}[\tilde{Z}_L(t)\leq c]\,dt.
\end{equation}
Define the modified Bessel function with index $\nu\in\mathbb{R}$ by 
\begin{equation} \label{def.besselfn}
I_\nu(x)=\left(\frac{x}{2}\right)^\nu \sum_{n=0}^\infty \frac{(x/2)^{2n}}{n!\Gamma(\nu+n+1)}, \quad x\geq 0.
\end{equation}
Clearly 
\begin{equation}\label{eq.besselfnasy}
\lim_{x\to 0^+} \frac{I_\nu(x)}{x^{\nu}} = \left(\frac{1}{2}\right)^\nu \frac{1}{\Gamma(\nu+1)}.
\end{equation}
Define $\delta=2\mu +1>0$ and $\nu(\delta)=\delta/2-1$. Since $\tilde{Z}_L$ obeys \eqref{eq.zltilde} we have 
\begin{equation}\label{eq.dfztildelt}
\mathbb{\tilde{P}}[\tilde{Z}_L(t)\leq c]=\int_0^c q_t(x_0^2,y)\,dy,
\end{equation}
where 
\begin{equation}\label{eq.densityztildelt} 
q_t(x_0^2,y)=\left\{
\begin{array}{cc}
 \frac{1}{2t}\left(\frac{y}{x_0^2}\right)^{\nu(\delta)/2}\exp\left(-\frac{x_0^2+y}{2t}\right)I_{\nu(\delta)}(|x_0|\sqrt{y}/t),  & x_0\neq 0,\\
 \frac{1}{(2t)^{\delta/2}}\frac{1}{\Gamma(\delta/2)} y^{\delta/2-1} \exp\left(-\frac{y}{2t}\right), & x_0=0.
\end{array}
\right.
\end{equation}
See e.g.,~\cite[Ch. XI, Corollary 1.4]{R&Y} and~\cite[Chap IV.,Example 8.3]{IkWat}, noting that there is a missing factor of $1/t$ in the first
formula in~\cite[Ch. XI, Corollary 1.4]{R&Y}. 

We now estimate $q_t(x_0^2,y)$, and use this to prove that $\mathbb{P}[\tilde{Z}_L(t)\leq c]$ tends to zero sufficiently 
quickly as $t\to\infty$. This guarantees that the integral on the righthand side of \eqref{eq.uprobzlprob} is finite and hence that 
\eqref{eq.uprobint} holds.   

In the case that $x_0=0$, since $\delta>0$ we have 
\begin{align*}
\mathbb{\tilde{P}}[\tilde{Z}_L(t)\leq c]
&=\frac{1}{(2t)^{\delta/2}}\frac{1}{\Gamma(\delta/2)} \int_0^c  y^{\delta/2-1} \exp\left(-\frac{y}{2t}\right)\,dy\\
&\leq \frac{1}{(2t)^{\delta/2}}\frac{1}{\Gamma(\delta/2)} \int_0^c  y^{\delta/2-1}\,dy 
= \frac{1}{(2t)^{\delta/2}}\frac{1}{\Gamma(\delta/2)} \frac{c^{\delta/2}}{\delta/2}. 
\end{align*}
Inserting this estimate into \eqref{eq.uprobzlprob} and using the fact that $\delta>0$ implies  \eqref{eq.uprobint}.

In the case when $x_0\neq 0$, by \eqref{eq.besselfnasy} there exists $\bar{x}_\delta>0$ such that 
\[
I_{\nu(\delta)}(x) < 2 \left(\frac{1}{2}\right)^{\nu(\delta)} \frac{1}{\Gamma(\nu(\delta)+1)}\cdot x^{\nu(\delta)}, \quad x<\bar{x}_\delta
\]  
Define $t_{\delta,c}=|x_0|\sqrt{c}/\bar{x}_\delta$. Then for $0\leq y\leq c$ and $t>t_{\delta,c}$ we have $|x_0|\sqrt{y}/t<\bar{x}_\delta$. 
Thus
\begin{equation*}
q_t(x_0^2,y)\leq \left(\frac{1}{2}\right)^{\delta/2-1}\frac{1}{\Gamma(\delta/2)}y^{\delta/2-1}\frac{1}{t^{\delta/2}}, \quad t>t_{\delta,c}, \,y\in[0,c].
\end{equation*} 
Therefore for $t\geq t_{\delta,c}$ the last estimate yields
\begin{align*}
\mathbb{\tilde{P}}[\tilde{Z}_L(t)\leq c]
&=\int_0^c q_t(x_0^2,y)\,dy
\leq \left(\frac{1}{2}\right)^{\delta/2-1}\frac{1}{\Gamma(\delta/2)} \frac{1}{t^{\delta/2}} 
\int_0^c y^{\delta/2-1}\,dy\\
&=\left(\frac{1}{2}\right)^{\delta/2-1}\frac{1}{\Gamma(\delta/2)} \frac{1}{t^{\delta/2}} \frac{c^{\delta/2}}{\delta/2}.
\end{align*}
Inserting this estimate into \eqref{eq.uprobzlprob} and using the fact that $\delta>0$ implies  \eqref{eq.uprobint}.

\subsection{Proof of Step B, i.e., \eqref{eq.conv1}}
We need to prove that if $U$ is the process defined in \eqref{eq.u}, then \eqref{eq.uprobint} implies
\eqref{eq.conv1}. Define $\pi^{\ast}_c(t)=\mathbb{\tilde{P}}[U(t)\leq ce^{-t}]$. Then by \eqref{eq.dfztildelt} we have 
\[
\pi^{\ast}_c(t)=\mathbb{\tilde{P}}[e^{-t}\tilde{Z}(e^{t}-1)\leq ce^{-t}]=\mathbb{\tilde{P}}[\tilde{Z}_L(e^{t-1})\leq c] = \int_0^c q_{e^t-1}(x_0^2,y)\,dy
\]
where $q$ is given by \eqref{eq.densityztildelt}. Therefore $\pi^\ast_c \in C([1,\infty),(0,\infty))$ and by \eqref{eq.uprobint} we have 
$\pi^\ast_c \in L^1([1,\infty),(0,\infty))$. Next define 
\[
F_c(t)=\frac{1}{e^t} \int_1^t e^s I_{\{U(s)\leq ce^{-s}\}}\,ds, \quad t\geq 1. 
\]
Then $\lim_{t\to\infty} F_c(t)=0$ a.s. implies \eqref{eq.conv1}. 

Since $\mathbb{\tilde{E}}[F_c(t)]=e^{-t}\int_1^t e^s \pi^\ast_c(s)\,ds$ and $\pi^\ast_c \in L^1([1,\infty),(0,\infty))$, we have that $\int_1^\infty \mathbb{\tilde{E}}[F_c(t)]\,dt < +\infty$.
Since $\pi^\ast_c \in C([1,\infty),(0,\infty))$, we have that $\mathbb{\tilde{E}}[F_c(\cdot)]\in C([1,\infty),(0,\infty))$. 
Then by \cite[Lemma 2.3]{App:2004a}, there exists a deterministic and increasing sequence $(a_n(c))_{n\geq 0}$ with $a_0=1$ and $a_n(c)\to\infty$
as $n\to\infty$ such that 
\begin{equation} \label{eq.fcest}
\sum_{n=0}^\infty \mathbb{\tilde{E}}[F_c(a_n(c))]<+\infty. 
\end{equation}
Next, define $G_c(n)=\int_{a_n(c)}^{a_{n+1}(c)} I_{U(s)\leq ce^{-s}}\,ds$. Then 
$\mathbb{\tilde{E}}[G_c(n)]=\int_{a_n(c)}^{a_{n+1}(c)} \pi^\ast_c(s)\,ds$. Since $\pi^\ast_c \in L^1([1,\infty),(0,\infty))$ 
we have 
\begin{equation} \label{eq.gcest}
\sum_{n=0}^\infty \mathbb{\tilde{E}}[G_c(n)]=\int_1^\infty \pi^\ast_c(s)\,ds<+\infty. 
\end{equation}
Now let $t\in[a_n(c),a_{n+1}(c)]$. Then
\begin{align*}
F_c(t)&=e^{-(t-a_n(c))}F_c(a_n(c))+\int_{a_n(c)}^t e^{-(t-s)} I_{\{U(s)\leq ce^{-s}\}}\,ds\\
&\leq F_c(a_n(c))+ G_c(n). 
\end{align*}
Therefore by \eqref{eq.fcest} and \eqref{eq.gcest} we have 
\[
\sum_{n=0}^\infty \mathbb{\tilde{E}}\left[ \sup_{a_n(c)\leq t\leq a_{n+1}(c)} F_c(t) \right] 
\leq 
\sum_{n=0}^\infty \mathbb{\tilde{E}}[F_c(a_n(c))]
+
\int_1^\infty \pi^\ast_c(s)\,ds < +\infty
\]
Therefore 
\[
\sum_{n=0}^\infty \sup_{a_n(c)\leq t\leq a_{n+1}(c)} F_c(t) < +\infty, \quad\text{a.s.}
\]
and so $\lim_{n\to\infty} \sup_{a_n(c)\leq t\leq a_{n+1}(c)} F_c(t)=0$ a.s. Hence $F_c(t)\to 0$ as $t\to\infty$ a.s. 

\subsection{Proof of Step C i.e., \eqref{eq.conv3}}
We now show that \eqref{eq.conv1} implies \eqref{eq.conv3}. By the definition of $\theta$ in \eqref{eq.theta}, we have that 
$\theta'(t)=g^2(\sqrt{Z_L(t)})$ for $t>0$. Recall moreover $0<K_1^2\leq g^2(x)\leq K_2^2$ for $x\in\mathbb{R}$. Therefore
\begin{align*}
\frac{1}{t}\int_1^t I_{\{Z_L(s)\leq c\}}\,ds
&=\frac{1}{t}\int_1^t I_{\{\tilde{Z}_L(\theta(s))\leq c\}}\,ds
=\frac{1}{t}\int_{\theta(1)}^{\theta(t)} I_{\{\tilde{Z}_L(u)\leq c\}}\tau'(u)\,du\\
&=\frac{\theta(t)}{t}\frac{1}{\theta(t)}\int_{\theta(1)}^{\theta(t)} I_{\{\tilde{Z}_L(u)\leq c\}}\frac{1}{g^2(\sqrt{\tilde{Z}_L(u)})}\,du\\
&\leq \frac{K_2^2}{K_1^2} \frac{1}{\theta(t)}\int_{\theta(1)}^{\theta(t)} I_{\{\tilde{Z}_L(u)\leq c\}}\,du,  
\end{align*}
where we used \eqref{eq.thetabounds} at the last step. Hence
\begin{equation} \label{eq.conv2a}
\limsup_{t\to\infty}
\frac{1}{t}\int_1^t I_{\{Z_L(s)\leq c\}}\,ds
\leq 
\frac{K_2^2}{K_1^2} \limsup_{t\to\infty} \frac{1}{t}\int_{0}^{t} I_{\{\tilde{Z}_L(u)\leq c\}}\,du. 
\end{equation}
Now using the connection between $\tilde{Z}_L$ and $U$ we have  
\[
 \frac{1}{t}\int_{0}^{t} I_{\{\tilde{Z}_L(u)\leq c\}}\,du=
  \frac{1}{t}\int_{0}^{\log(1+t)} I_{\{U(s)\leq ce^{-s}\}} e^s\,ds. 
\]
Using this identity and \eqref{eq.conv2a} we have 
\begin{align*}
\limsup_{t\to\infty}
\frac{1}{t}\int_1^t I_{\{Z_L(s)\leq c\}}\,ds
&\leq
\frac{K_2^2}{K_1^2} 
\limsup_{t\to\infty} \frac{1}{t}\int_{0}^{\log(1+t)} I_{\{U(s)\leq ce^{-s}\}} e^s\,ds\\
&=
\frac{K_2^2}{K_1^2}
\limsup_{t\to\infty} \frac{1}{e^t}\int_{1}^{t} I_{\{U(s)\leq ce^{-s}\}} e^s\,ds.
\end{align*}
Therefore, as \eqref{eq.conv1} holds, we have \eqref{eq.conv3}. 

\subsection{Proof of Step D i.e., \eqref{eq.thetaaverage}}
We have that \eqref{eq.conv3} holds i.e., for each $c>0$ there exists an a.s. event $\Omega_c$ such that  
\[
\Omega_c:=\left\{\omega: 
\lim_{t\to\infty} \frac{1}{t}\int_0^t I_{\{Z_L(s,\omega)\leq c\}}\,ds = 0\right\}. 
\]
Moreover, for $\omega \in \Omega_c$ we also have 
\begin{equation}\label{eq.indicatoron}
\lim_{t\to\infty} \frac{1}{t}\int_0^t I_{\{Z_L(s,\omega)> c\}}\,ds = 1. 
\end{equation}

For every $\varepsilon\in(0,1)$, there exists $c(\varepsilon)>0$ such that $\sigma^2(1-\varepsilon)< g^2(x) < \sigma^2(1+\varepsilon)$   
for all $x>\sqrt{c(\varepsilon)}$. Therefore
\[
\sigma^2(1-\varepsilon)I_{\{Z_L(s)>c(\varepsilon)\}}\leq g^2(\sqrt{Z_L(s)})I_{\{Z_L(s)>c(\varepsilon)\}}
\leq
\sigma^2(1+\varepsilon)I_{\{Z_L(s)>c(\varepsilon)\}}.
\]
Therefore 
\begin{multline}\label{eq.indicatorupperlower} 
\sigma^2(1-\varepsilon)\frac{1}{t}\int_0^t I_{\{Z_L(s)>c(\varepsilon)\}}\,ds \leq \frac{1}{t}\int_0^t g^2(\sqrt{Z_L(s)})I_{\{Z_L(s)>c(\varepsilon)\}}\,ds
\\
\leq
\sigma^2(1+\varepsilon) \frac{1}{t}\int_0^t I_{\{Z_L(s)>c(\varepsilon)\}}\,ds.
\end{multline}
Now, for $\omega\in \Omega_{c(\varepsilon)}$, by using \eqref{eq.indicatoron} and the lefthand member of \eqref{eq.indicatorupperlower}, 
we have 
\begin{align*}
\liminf_{t\to\infty} \frac{1}{t}\int_0^t g^2(\sqrt{Z_L(s)})\,ds
&\geq 
\liminf_{t\to\infty} \frac{1}{t}\int_0^t g^2(\sqrt{Z_L(s)})I_{\{Z_L(s)>c(\varepsilon)\}}\,ds\\
&\geq \sigma^2(1-\varepsilon)
\liminf_{t\to\infty} \frac{1}{t}\int_0^t I_{\{Z_L(s)>c(\varepsilon)\}}\,ds\\
&=\sigma^2(1-\varepsilon).
\end{align*}
Therefore with $\Omega_1^\ast=\cap_{\varepsilon\in(0,1)\cap\mathbb{Q}}\Omega_{c(\varepsilon)}$ we have that $\Omega_1^\ast$ is an almost sure 
event and moreover
\begin{equation} \label{eq.lowerthetaaverage}
\liminf_{t\to\infty} \frac{1}{t}\int_0^t g^2(\sqrt{Z_L(s)})\,ds\geq \sigma^2, \quad\text{a.s. on $\Omega_1^\ast$}.
\end{equation}

To obtain an upper bound, first note that $g^2(x)\leq K_2^2$ implies 
\[
\frac{1}{t}\int_0^t g^2(\sqrt{Z_L(s)})I_{\{Z_L(s)\leq c(\varepsilon)\}}\,ds
\leq K_2^2
\frac{1}{t}\int_0^t I_{\{Z_L(s)\leq c(\varepsilon)\}}\,ds. 
\]
Therefore 
\begin{equation} \label{eq.upperaveragethetaest1}
\limsup_{t\to\infty}
\frac{1}{t}\int_0^t g^2(\sqrt{Z_L(s)})I_{\{Z_L(s)\leq c(\varepsilon)\}}\,ds=0, \quad \text{a.s. on $\Omega_{c(\varepsilon)}$}.
\end{equation}
Since 
\begin{multline*}
\frac{1}{t}\int_0^t g^2(\sqrt{Z_L(s)})\,ds
\\=
 \frac{1}{t}\int_0^t g^2(\sqrt{Z_L(s)})I_{\{Z_L(s)\leq c(\varepsilon)\}}\,ds
+
 \frac{1}{t}\int_0^t g^2(\sqrt{Z_L(s)})I_{\{Z_L(s)> c(\varepsilon)\}}\,ds,
\end{multline*}
by combining \eqref{eq.upperaveragethetaest1} with the righthand member of \eqref{eq.indicatorupperlower} we get 
\[
\limsup_{t\to\infty}
\frac{1}{t}\int_0^t g^2(\sqrt{Z_L(s)})\,ds
\leq \sigma^2(1+\varepsilon) \quad \text{a.s. on $\Omega_{c(\varepsilon)}$}.
\] 
Therefore with $\Omega_1^\ast$ as the almost sure event defined above we have 
\begin{equation} \label{eq.upperthetaaverage}
\liminf_{t\to\infty} \frac{1}{t}\int_0^t g^2(\sqrt{Z_L(s)})\,ds\leq \sigma^2, \quad\text{a.s. on $\Omega_1^\ast$}.
\end{equation}
Combining \eqref{eq.lowerthetaaverage} and \eqref{eq.upperthetaaverage} yields \eqref{eq.thetaaverage} as required.


\end{document}